\documentclass[alpha-refs]{wiley-article}

\usepackage{graphicx}
\usepackage[space]{grffile}
\usepackage{latexsym}
\usepackage{textcomp}
\usepackage{longtable}
\usepackage{tabulary}
\usepackage{booktabs,array,multirow}
\usepackage{amsfonts,amsmath,amssymb}
\usepackage{natbib}
\usepackage{url}
\usepackage{hyperref}
\hypersetup{colorlinks=false,pdfborder={0 0 0}}
\usepackage{etoolbox}
\usepackage{algpseudocode}
\usepackage{algorithm}
\usepackage{xcolor}
\usepackage{comment}

\newif\iflatexml\latexmlfalse

\AtBeginDocument{\DeclareGraphicsExtensions{.pdf,.PDF,.eps,.EPS,.png,.PNG,.tif,.TIF,.jpg,.JPG,.jpeg,.JPEG}}

\usepackage[utf8]{inputenc}
\usepackage[english]{babel}

\usepackage{siunitx}

\iflatexml

\else

\corraddress{Koki Ho PhD, Daniel Guggenheim School of Aerospace Engineering, Georgia Institute of Technology, Atlanta, Georgia, 30068, United States}
\corremail{kokiho@gatech.edu}
\presentadd{Daniel Guggenheim School of Aerospace Engineering, Georgia Institute of Technology, Atlanta, Georgia, 30068, United States}
\fundinginfo{National Science Foundation, Grant/Award Number: 1942559}
\fi

\papertype{Original Article}

\title{Optimizing Flexible Complex Systems with Coupled and Co-Evolving Subsystems under Operational Uncertainties}

\author[1]{Koki Ho}
\author[1]{Masafumi Isaji}
\author[1]{Malav Patel}
\author[1]{Kayla Garoust}
\affil[1]{Georgia Institute of Technology}

\runningauthor{Koki Ho}

\begin{document}

\maketitle
\selectlanguage{english}
\begin{abstract}
The paper develops a novel design optimization framework and associated computational techniques for staged-deployment optimization of complex systems under operational uncertainties. It proposes a local scenario discretization method that offers a computationally efficient approach to optimize staged co-deployment of multiple coupled subsystems by decoupling weak dynamic interaction among subsystems. With case studies, the proposed method is demonstrated to provide an effective and scalable optimization strategy to determine the flexible systems design under uncertainty. Extensions of the proposed method in the context of system partitioning are also explored. The developed optimization framework is expected to improve the staged-deployment design of various complex engineering systems, such as water, energy, food, and other infrastructure systems.

\end{abstract}%

\section{Introduction}
With dynamic variations in the future operating environment and the rapid evolution of technology, traditional approaches to designing and deploying a complex engineering system all at once can involve significant financial and technical risks. Instead, deploying the system progressively in stages, starting from affordable stages, has been shown to be an effective solution, as the system can adapt to operational uncertainties over time and thus mitigate risks \cite{neufville2011flexibility,de2004staged,hamdan2023dynamic}. A successful example of a staged-deployed system is the Blue Cross Blue Shield Tower in Chicago, which completed its first stage as a 33-story building in 1997 with built-in flexibility (e.g., extra steel and elevator shafts) for expansion, and later completed its second stage as a 57-story building in 2010 in response to the increasing demand for office space.

Conventionally, the flexibility gained by incorporating evolutionary design options has been analyzed using real options through tradespace exploration or stochastic programming \cite{de2004staged,hassan2005value,engel2008designing}. With few exceptions where closed-form solutions are available, general flexible design problems have been tackled by discretizing the probability distributions of the time-variant operational uncertainties into scenarios and analyzing the impact of flexible decisions over each scenario. In the mathematical programming field, sample average approximation methods \cite{kleywegt2002sample,kim2014guide,homem2014monte} are often used, and various scenario generation and reduction methods exist using probability metrics \cite{dupavcova2003scenario,heitsch2003scenario,rujeerapaiboon2022scenario,hoyland2003heuristic} or problem structures \cite{henrion2022problem,bertsimas2023optimization,hewitt2022decision}. In the design science literature, binomial lattice is also often used to generate the scenarios \cite{neufville2011flexibility}. Given the discretized scenarios, stochastic programming with a standard formulation or a decision-rule-based formulation \cite{caunhye2017approach,cardin2017approach} is used to optimize the flexible design decisions. 

However, for a complex system composed of coupled subsystems, optimizing such flexible system designs is extremely challenging because we need effective coordination between the design and deployment of each subsystem in response to uncertainties. Consider an application of designing and deploying an infrastructure system consisting of a water generation system and an electric power plant in a developing area. These two infrastructure elements (i.e., ``subsystems'') are interdependent because operating the water generation system requires electricity, and operating the power plant requires
water \cite{siddiqi2011water}. Therefore, we need a co-deployment optimization for these infrastructure elements, considering not only their own demand uncertainties but also their interdependence. This problem becomes even more complex when more coupled subsystems, such as a food production system, and their corresponding uncertainties are introduced. When there are many such interdependent subsystems to be co-deployed, we do not currently have a computationally scalable method that can co-optimize these subsystem-level designs and staged-deployment decisions in response to their operational uncertainties. With the traditional methods, we would need to discretize the joint probability distribution of all subsystems’ uncertainties into a (combinatorially large) set of scenarios and analyze those scenarios for all decisions in all subsystems in the optimization, resulting in an extremely large number of variables and constraints. This would make the optimization poorly scalable and computationally challenging. Existing literature circumvented this challenge by designing the stages with a top-down approach using simple system-level low-fidelity models or experts’ opinions, which often lack scalability and rigor. Thus, there is an imperative need for a bottom-up coordinated subsystem-level design method with a scalable optimization formulation while providing a sufficiently flexible design solution. 

This paper addresses this imperative need for flexible, complex multi-subsystem system design by creating a unique design optimization framework of coordinated staged development. The proposed design optimization framework examines the staged-deployment strategies as well as the needed ``built-in flexibility'' to enable later flexible deployment in response to multiple sources of uncertainty. Our key idea is to decouple weak dynamic subsystem interaction so that the built-in flexibility can be designed at the subsystem level, while considering only the strong couplings from other subsystems under uncertainties. This approach is enabled by our developed method to quantify and decouple the weak dynamic subsystem couplings. The developed method is demonstrated with case studies of two-stage multi-subsystem staged deployment to illustrate its computational cost savings while maintaining the flexibility in the design solution. Furthermore, extensions of the proposed method in the context of system partitioning are also investigated.

The remainder of this paper is structured as follows. Section \ref{Formulation1} formulates the conventional stochastic programming formulation for flexible design of systems, Section \ref{Formulation2} extends that formulation for complex multi-subsystem systems, and Section \ref{method} discusses the newly proposed local scenario discretization method. Furthermore, Section \ref{demo} provides the demonstration case studies, followed by Section \ref{discussion} that discusses the implications and extensions of the proposed method, particularly in the context of large-scale multi-subsystem problems. Finally, Section \ref{conclusion} summarizes the main contribution of the paper.

\section{Overview of Staged-Deployment Problem Formulation}
\label{Formulation1}
Conventionally, the flexible staged-deployment design for engineering systems can be formulated as a multi-stage nonlinear stochastic programming problem. We use scenarios to represent uncertainties at each discrete decision time step. Note that for deployment decisions under uncertainty, discretization of decision time steps is a reasonable and common assumption \cite{ahmed2003MultiStageStochasticInteger}. At each discrete time step, we also assume that the support of random variables is finite; that is, the number of considered scenarios is finite. This finite support assumption is also a standard practice \cite{ahmed2004TSSIP_BnB,vanderlaan2023MIRbenders}. In case a random variable follows a continuous distribution, it can be discretized into a finite set of scenarios. For a multi-stage problem, earlier-stage decisions have less information about the uncertainties than later-stage ones; this concept is commonly known as nonanticipativity in Stochastic Programming \cite{shapiro2021lec_on_SP}. 
For example, the Stage 1 decision does not depend on the scenario because no information about the scenario is available then, while the Stage 2 decision would respond to the specific scenario realized over time. This paper focuses on such a two-stage staged-deployment problem.

Mathematically, a two-stage staged-deployment problem can be formulated as follows: denote the domain of $x$ by $\mathcal{X}$, and denote the scenario set by $\mathcal{S}$, defined by discretizing the variation of parameters $D$. In this paper, we assume that all scenarios have the same probability and are ordered in ascending $D$ value, but the proposed method can be extended to a more general case. Each variable $x$ belongs to a decision for Stage 1, $x^1$  (independent of scenario $s$), or for Stage 2, $x^{2s}$ (dependent on scenario $s$), where the stage and scenario are represented in superscript. In addition, we also define the first and second-stage cost functions $C^1$ and $C^{2s}$; here, $C^1$ is dependent on the decision at Stage 1, $x^1$, where $C^{2s}$ is dependent on both the decision at Stages 1 and 2 in its corresponding scenario realization, $x^1, x^{2s}$. The general two-stage formulation can be written as:
\begin{gather}
    \label{eq:single}
    \begin{aligned}
        \min_{x\in\mathcal{X}} \quad &C
        \\
        \text{s.t.}\quad 
        \\
        C&=C^1 (x^1)+\mathrm{E}_{s\in\mathcal{S}} [C^{2s} (x^1,x^{2s})]
    \end{aligned}
\end{gather}

As an example, consider a capacity design problem for an infrastructure system (e.g., electricity power plant) with the following formulation. 
\begin{gather}
    \label{eq:single_example}
    \begin{aligned}
        \min_{x\in\mathcal{X}} \quad &C
        \\
        \text{s.t.}\quad 
        \\
        C&=c^1 \cdot (x^1)^\alpha+\mathrm{E}_{s\in\mathcal{S}} [c^2 \cdot (dx^{2s})^\alpha] \\
        dx^{2s}&=\max\{x^{2s}-x^1,0\}\quad\forall s\in\mathcal{S}\\
        x^1&\geq D^1\\
        x^{2s}&\geq D^{2s}\quad\forall s\in\mathcal{S}
        \end{aligned}
\end{gather}
In this context, $D$ is a positive parameter representing the demand of the system; specifically, we know the Stage 1 demand $D^1$, but the Stage 2 demand $D^{2s}$ is dependent on scenario realization $s$ and therefore still uncertain. Our decision variable $x^1$ corresponds to the initial infrastructure capacity at Stage 1, which should be at least more than sufficient to satisfy the Stage 1 demand, but can also have additional capacity to anticipate future demand growth. On the other hand, $x^{2s}$ corresponds to the capacity decision at Stage 2. This decision needs to satisfy the Stage 2 demand, which may be more or less than the Stage 1 demand. If Stage 2 demand is less than the existing capacity from Stage 1, $x^1$, then no action is needed at Stage 2 because the system can handle the demand at Stage 2 without expansion. However, if the Stage 1 capacity $x^1$ is not sufficient to satisfy the demand at Stage 2 under scenario $s$, the capacity must be expanded by $x^{2s}- x^1$. As such, the required expansion size can therefore be expressed as $dx^{2s}=\max\{x^{2s}-x^1,0\}$. For the considered optimization problem, this equality constraint involving the $\max$ function can be substituted with an inequality one, $dx^{2s}\geq \max\{x^{2s}-x^1,0\}$, which can be expressed as a set of linear inequalities. The cost function in this example is assumed to be a parametric function, $C^1(x^1)=c^1 \cdot (x^1)^\alpha$ and $C^{2s}(x^1,x^{2s})=c^2 \cdot (dx^{2s})^\alpha=c^2 \cdot  (\max\{x^{2s}-x^1,0\})^\alpha$ where $c^1, c^2 > 0$ are predetermined deployment cost coefficients for Stage 1 and 2, respectively. Because deploying later costs more than deploying earlier (i.e., $c^1 < c^2$), we typically add some excess capacity to the Stage 1 design beyond its demand (i.e., built-in flexibility). However, since the Stage 2 demand can be larger or smaller than the Stage 1 demand, the optimizer performs a tradeoff between adding more capacity at Stage 1 to prepare for a potentially growing demand vs. keeping the Stage 1 capacity limited to prepare for a potentially flat or declining demand and minimize the potential waste of the unused capacity. This tradeoff is performed by minimizing the expected cost over the possible scenarios under the constraints that the capacity is equal to or larger than the demand for both stages and for all scenarios.

\section{Staged Co-Deployment of Multi-Subsystem Systems}
\label{Formulation2}
The conventional formulation introduced in the previous section is now extended to more complex systems with multiple coupling subsystems. Specifically, we consider a set of subsystems $\mathcal{Z}$, and each Subsystem $i \in \mathcal{Z}$ has its own uncertain parameter $D_i$ and decision variables $x_i$ for each stage. The subsystems are dependent on each other through their coupling variables. Namely, the variables now include the local variable $x_i$ for each Subsystem $i$ and the coupling variables $x_{ij}$ from subsystems $j$ to $i$. Similarly to the single-subsystem case, each local and coupling variable belongs to a decision for Stage 1, $x^1$ (independent of scenario $s$), or for Stage 2, $x^{2s}$ (dependent on scenario $s$). The general two-stage formulation can be written as:
\begin{gather}
    \label{eq:full}
    \begin{aligned}
        \min_{x\in\mathcal{X}}\quad &\sum_{i\in \mathcal{Z}} C_i
        \\
        \text{s.t.}\quad &
        \\
        C_i&=C_i^1 (x_i^1,x_{ij:j\neq i}^1)+\mathrm{E}_{s\in\mathcal{S}} [C_i^{2s} (x_i^1,x_{ij:j\neq i}^1,x_i^{2s},x_{ij:j\neq i}^{2s})] 
    \end{aligned}
\end{gather}
The difference between this formulation and the previous one is that it considers the coupling between different subsystems. The cost function for each subsystem is also dependent on both its own local variables and the coupling variables from all other subsystems. 

Consider again the example of infrastructure capacity expansion, but this time, for two subsystems (e.g., water generation infrastructure and electricity power plant infrastructure), denoted by Subsystems A and B. We assume that there is a coupling between these two subsystems; for example, a larger water generation infrastructure capacity leads to a larger electricity demand, and a larger electricity power plant leads to a larger water demand. These couplings are captured by the coupling variables $x_{AB}$ (from B to A) and $x_{BA}$ (from A to B) as shown in Fig. \ref{coupling}.

The formulation is shown below:
\begin{gather}
    \label{eq:full_example}    
    \begin{aligned}
        \min_{x\in\mathcal{X}}\quad &C_A+C_B
        \\
        \text{s.t.}\quad &
        \\
        \text{Subsystem A: } &\text{(Receives $x^1_{AB},x^{2s}_{AB}$ from Subsystem B)}\\        
        C_A&=c^1_A \cdot (x^1_A)^\alpha+\mathrm{E}_{s\in\mathcal{S}} [c^{2}_A \cdot (dx^{2s}_A)^\alpha] \\
        dx^{2s}_A&=\max\{x^{2s}_A-x^1_A,0\} \quad\forall s\in\mathcal{S}\\
        x^1_A&\geq D^1_A+d_{AB} \cdot x^1_{AB}\\
        x^{2s}_A&\geq D^{2s}_A+d_{AB} \cdot x^{2s}_{AB} \quad\forall s\in\mathcal{S}\\
        x^1_{BA}&=x^1_{A}\\
        x^{2s}_{BA}&=x^{2s}_A \quad\forall s\in\mathcal{S}\\
        \text{Subsystem B: } &\text{(Receives $x^1_{BA},x^{2s}_{BA}$ from Subsystem A)}\\        
        C_B&=c^1_B \cdot (x^1_B)^\alpha+\mathrm{E}_{s\in\mathcal{S}} [c^2_B \cdot (dx^{2s}_B)^\alpha] \\
        dx^{2s}_B&=\max\{x^{2s}_B-x^1_B,0\} \quad\forall s\in\mathcal{S}\\
        x^1_B&\geq D^1_B+d_{BA} \cdot x^1_{BA}\\
        x^{2s}_B&\geq D^{2s}_B+d_{BA} \cdot x^{2s}_{BA} \quad\forall s\in\mathcal{S}\\
        x^1_{AB}&=x^1_{B}\\
        x^{2s}_{AB}&=x^{2s}_B \quad\forall s\in\mathcal{S}\\
    \end{aligned}
\end{gather}

\begin{figure}[h!]
\begin{center}
\includegraphics[width=0.50\columnwidth]{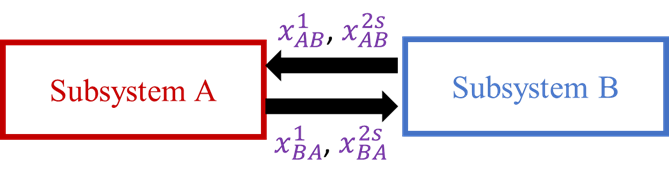}
\caption{{Coupling variables between the staged-deployed subsystems.
{\label{coupling}}%
}}
\end{center}
\end{figure}

This formulation captures the co-deployment of both subsystems A and B under uncertainty in both demands. Specifically, for both subsystems, Stage 1 demands, $D_A^1$ and $D_B^1$, are known, but Stage 2 demands, $D_A^{2s}$ and $D_B^{2s}$, are uncertain, represented by scenarios $s\in\mathcal{S}$. If we consider $S_A$ scenarios for Subsystem A’s demands and $S_B$ scenarios for Subsystem B’s demands, all the $S_A \cdot  S_B$ combinations will be included in the scenario set $\mathcal{S}$ due to the coupling between these two subsystems. In the same way as the single-subsystem case, the variables are defined as the capacity for each subsystem for the first stage, $x_A^1$ and $x_B^1$, and for the second stage in each scenario, $x_A^{2s}$ and $x_B^{2s}$. Using the same cost model, deploying Stage 1 with the size of $x_A^1$ and $x_B^1$ requires cost $c_A^1 \cdot (x_A^1)^\alpha$ and $c_B^1 \cdot (x_B^1)^\alpha$, respectively, and expanding its capacity by $dx_A^{2s}$ and $dx_B^{2s}$ also requires additional cost $c_A^2 \cdot (dx_A^{2s})^\alpha$ and $c_B^2 \cdot (dx_B^{2s})^\alpha$, respectively, when necessary; when the needed capacity for Stage 2 is the same as or smaller than the Stage 1 capacity, we do not need extra expansion cost. The demand for each subsystem includes both the demand of each subsystem itself (e.g., $D_A^1$ and $D_A^{2s}$) and the demand caused by the other subsystem (e.g., the terms $d_{AB} \cdot x_{AB}^1$ and $d_{AB} \cdot x_{AB}^{2s}$ represent the coupling effect from the origin subsystem B to the destination subsystem A, where $d_{AB}$ is a positive parameter). For example, any capacity decision in Subsystem B, $x_B$, impacts Subsystem A through a coupling variable $x_{AB}$ and adds its demand by $d_{AB}\cdot x_{AB}$, and vice versa, for both Stages and all scenarios. Due to this coupling, we need to design the co-deployment decision for both subsystems, considering both sources of uncertainty in these two subsystems. This problem can be solved either as a single monolithic optimization problem or through problem decomposition into smaller problems solved in parallel and iteratively. This formulation is referred to as the \textit{fully flexible} formulation in this paper. 

 Note that this paper assumes the source of uncertainty in each subsystem is independent, but the flexible decisions responding to each subsystem's source of uncertainty can impact the decisions in other subsystems due to the inter-subsystem dependency. One way to interpret this assumption is that we break the uncertainty in each subsystem down to an independent component (e.g., uncertain electricity demand) and a coupled component (e.g., additional electricity demand caused by the increased water infrastructure capacity responding to uncertain water demand), the latter of which is modeled as a coupling between subsystems. This is in contrast to other methods that directly leverage the dependency between the uncertainties, such as variance reduction in Monte Carlo methods \cite{homem2014monte,homem2014stochastic}.

While this fully flexible design formulation can be extended to more subsystems and more sources of uncertainty, it is evident that the number of scenarios that need to be considered will increase exponentially. In this formulation, we would need to consider all combinations of uncertainty scenarios for all subsystems to optimize the flexible decisions, which can become computationally prohibitive. This poor scalability can be a significant challenge for large-scale complex systems design and deployment when there are multiple sources of uncertainty. In this paper, we address this challenge by developing a more efficient \textit{local} scenario discretization method, as detailed in the next section, and improve the computational efficiency.

\section{Local Scenario Discretization Method}
\label{method}
\subsection{Overview and Related Works}
The solution method proposed in this paper addresses the scalability challenges of the fully flexible formulation in Eq. \eqref{eq:full}, which is due to the large set of scenarios as a result of the combinatorial enumeration of multiple independent uncertainties. In the conventional fully flexible approach, as reviewed previously, these scenarios are assigned and analyzed uniformly for all decision variables in all subsystems, resulting in a large number of variables and constraints. Despite the development in stochastic programming methods \cite{pflug2014multistage,homem2014monte}, such a formulation is intractable for large-scale systems. This poor scalability in the number of variables can be partially alleviated with a decision-rule-based formulation \cite{cardin2017approach}, but it requires prior knowledge of the form of the decision rule to be used, which can be difficult to obtain for complex systems. Approximate dynamic programming techniques have been proposed to alleviate the curse of dimensionality in the time dimension if the problem can be formulated as a dynamic programming problem \cite{powell2007approximate}, but complexity due to the number of uncertainty sources remains a computational challenge. Scenario reduction methods have been developed to reduce the number of scenarios via probability distance or moment-matching \cite{dupavcova2003scenario,heitsch2003scenario,rujeerapaiboon2022scenario,hoyland2003heuristic}, but they do not consider the relative significance among the sources of uncertainties and therefore are not suitable for the design of complex systems with a variety of uncertainty sources. More recent works focused on scenario reduction using cost functions and problem structures \cite{henrion2022problem,bertsimas2023optimization,hewitt2022decision}, but these works typically need to solve a separate optimization problem to define the scenarios, which might not be computationally feasible in the context of complex systems design. Furthermore, surrogate- or machine learning-assisted scenario generation and reduction methods have also been developed \cite{wu2022learning,patel2022neur2sp,suemitsu2025surrogate}, but they require surrogate model training, which, again, might not be realistic for large systems. In practice, system designers often filter out insignificant sources of uncertainty and focus only on the major ones; since this approach would fully include or fully exclude certain sources of uncertainty across the entire problem, the resulting solution would be limited in its flexibility as it can only consider a small subset of the uncertainty sources for a given computational limit. Thus, there is a need to develop a computationally efficient method that can manage the coordination complexity of a set of interacting uncertain subsystems with minimal sacrifice in design flexibility. 

In this paper, we develop a new local scenario discretization method for computationally efficient co-deployment management for complex systems. In this method, we partially or completely decouple weak uncertain couplings between subsystems to reduce computational cost, while largely maintaining design flexibility. Namely, instead of fully including or excluding each source of uncertainty across all variables as the conventional fully flexible approach does, we locally assign the scenarios to each subsystem's variables taking into account the strength of the couplings from other subsystems.

To illustrate the motivation, we revisit the example of the formulation in Eq. \eqref{eq:full_example}. In this problem, if we completely decouple a coupling variable $x_{BA}^{2s}$, the coupling from Subsystem A to B is written as a single variable $x_{BA}^2$ instead of as flexible variables that vary for each scenario $s$. 
With this decoupling, Subsystem B only needs to evaluate $S_B$ scenarios instead of $S_A\cdot S_B$ ones. Similarly, if we completely decouple both $x_{BA}^{2s}$ and $x_{AB}^{2s}$, then Subsystem A evaluates $S_A$ scenarios and Subsystem B evaluates $S_B$ scenarios, instead of $S_A \cdot S_B$ ones for both. While these examples assume complete decoupling of the weak coupling, we can further consider partial decoupling, i.e., assigning the number of scenarios for the coupling variables based on the strength of the coupling; this approach enables us to effectively maintain the flexibility of the design solution with a lower computational cost.

More generally, we develop a method to assign the number of scenarios to each coupling variable $x_{ij}^s$ (i.e., from subsystems $j$ to $i$) based on the impact of its uncertainty on the output of its destination Subsystem $z_i^s$, compared to the impact of Subsystem $i$'s own uncertain parameter(s) $D_i$; this metric is referred to as dynamic coupling $\delta_{ij}$ in this paper. The rationale behind this idea is to reduce the computational cost by considering only the variation of the coupling variables that significantly impact the destination subsystem’s output $z_i^s$. Note that this idea shares the guiding principle to some extent with the conventional approach that screens out insignificant uncertainty sources; however, our idea is novel in that it eliminates (or at least reduces) the impacts of insignificant uncertain couplings between the subsystems, instead of eliminating insignificant uncertainty sources themselves across the entire problem. This innovation enables us to efficiently allocate computational resources to the stronger couplings and thus is expected to improve the computational efficiency as compared to the conventional fully flexible approach that considers the same set of scenarios for all subsystems. 

While the above explanation provides a high-level overview of the proposed method, the steps to determine the dynamic coupling strength and the formulation of the locally discretized formulation are discussed in more detail in the following subsections.
 
\subsection{Dynamic Coupling Characterization via Local Sensitivity Analysis}
A key step in the proposed approach is to determine the number of scenarios for the coupling variables for the main optimization. In this step, we characterize the dynamic coupling $\delta_{ij}$, corresponding to the impact of the uncertainty of the coupling variable(s) $x_{ij}$ on the output $z_i$, relative to that of Subsystem $i$'s own uncertain parameter(s) $D_i$. Note that the uncertain input to a model can impact the total cost directly and indirectly. For direct impacts (i.e., main effects), the output corresponds to the cost function of the subsystem, whereas for indirect impacts, it corresponds to the coupling variables originating from that subsystem, which affect the total cost through other subsystems. As a first-order analysis, we determine dynamic coupling based only on the direct impacts.

While there are various methods to evaluate such impacts, in this paper, we propose an approach to determine the coupling strengths using the analytic expression of the sensitivity information. 
For example, the impact of the uncertain parameter within the Subsystem $i$, $D_i$, on $C_i$ can be represented by $\partial C_i/\partial D_{i}$, and the impact of the coupling variables $x_{ij}$ flowing into Subsystem $i$ on Subsystem $i$'s cost $C_i$ can be represented by $\partial C_i/\partial x_{ij}$.  
We assume that these quantities can be defined and an analytic relationship between these quantities that is common for all scenarios $s$ can be obtained. This information is used to determine how to discretize the scenarios for the coupling variables.

Specifically, we quantify the effect of a given coupling (i.e., the extent to which the coupling can be decoupled or reduced without significant impact on the accuracy of the cost evaluation) by comparing the (partial) derivatives of the cost with respect to the coupling variables coming from another subsystem, $\partial C_i/\partial x_{ij}$, and that with respect to the subsystem's own uncertainty source(s), $\partial C_i/\partial D_{i}$. For example, if $\partial C_i/\partial x_{ij}$ is half of $\partial C_i/\partial D_{i}$, one would need half as many scenarios over the same (absolute) value range for $x_{ij}$ as for $D^s_{i}$ to reach the same granularity in the output discretization. To capture a different scale of $x_{ij}$ and $D_{i}$, we need to normalize the derivative. Since we do not have the actual values of the variables until we solve the problem and we also need a common expression across all scenarios, we use the deterministic result (i.e., inflexible design for the conservative, worst-case demand) for normalization.  
Thus, we define the dynamic coupling as follows:
\begin{equation}
\delta_{ij}=\left|\frac{\partial C_i/\partial x_{ij}}{\partial C_i/\partial D_{i}}\cdot\frac{x_{ij}^\text{det}}{D_i^\text{det}}\right|
\end{equation}
in case both derivatives are nonzero, where the superscript ``det'' denotes the deterministic solution. Thus, since the subsystem's own uncertain parameter is discretized into $S_i$ scenarios, the minimum number of scenarios needed for a coupling variable $x_{ij}$ to capture its coupling impact can be written as:
\begin{equation}
S_{ij}=\lceil S_i \cdot \delta_{ij} \rceil
\end{equation}
Note that one interpretation of this dynamic coupling is the relative cost sensitivity to external vs. internal uncertainty, which captures the strength of the coupling between the subsystem of interest and other subsystems.

With the logic described above, the number of scenarios for each subsystem's decision can be defined in the following way. Considering all scenario combinations of its own uncertainty $D_i$ and the coupling uncertainty $x_{ij}$, the number of scenarios Subsystem $i$'s decision variables need to consider is:
\begin{equation}
    \sigma_i=S_i \cdot \prod_{j:j\neq i} S_{ij}
\end{equation}
To prevent the (rare) case where the number of scenarios becomes larger than in the fully flexible case (i.e., $\prod_{j} S_{j}$), we refine the definition of $S_{ij}$ as:
\begin{equation}
S_{ij}=\min\left\{\lceil S_i \cdot \delta_{ij} \rceil,S_i\right\}
\end{equation}
In the example formulation shown in Eq. \eqref{eq:full_example}, assuming $d_{ij}$, $x_{ij}^\text{det}$, and $D_i^\text{det}$ are positive,
\begin{equation}
\delta_{ij}=d_{ij}\cdot(x_{ij}^\text{det}/D_i^\text{det})
\end{equation}
at the optimal solution where the linear independence constraint qualification is satisfied and the Stage 2 demand constraints are active. Thus, we set 
\begin{equation}
    S_{ij}=\min\left\{\lceil S_i \cdot d_{ij}\cdot(x_{ij}^\text{det}/D_i^\text{det}) \rceil,S_i\right\}
\end{equation}
and therefore, the number of scenarios for the decision variables in Subsystem $i$ is 
\begin{equation}
    \sigma_i=S_i \cdot \prod_{j:j\neq i} \min\left\{\lceil S_i \cdot d_{ij}\cdot(x_{ij}^\text{det}/D_i^\text{det}) \rceil,S_i\right\}
\end{equation}

The number of scenarios for the coupling variable $x_{ij}$ in the constraints for Subsystem $i$ can be $S_{ij}$ (the minimum needed to capture the coupling impact) or $\sigma_i$ (the number used by the local decision variables $x_i$). When solving the entire optimization as a single monolithic problem, considering fewer than $\sigma_i$ scenarios for $x_{ij}$ does not necessarily reduce the size of the problem further because Subsystem $i$ needs to evaluate and consider the decisions for $\sigma_i$ scenarios anyway. Thus, this paper assumes $\sigma_i$ scenarios for the coupling variables $x_{ij}$.\footnote{In this case, the coupling variables would never be ``completely'' decoupled because the number of scenarios for the coupling variables is equal to the number of scenarios for the local variables, which is always larger than one unless the problem is deterministic.} However, when solving this problem using a decomposition-based approach, where each subsystem is solved in parallel and iteratively, we may prefer to minimize the number of scenarios for the coupling variables $x_{ij}$, in which case assigning $S_{ij}$ scenarios to $x_{ij}$ may be preferred at the expense of lost information due to further scenario reduction. 

A great benefit of this method is that the needed information can be obtained directly from the formulation, and therefore there is no need to run multiple random cases to obtain reliable results, unlike other alternative methods such as sample average approximation \cite{kleywegt2002sample}.

\subsection{Formulation of Locally Discretized Main Optimization Problem}
The proposed approach characterizes the dynamic strength $\delta_{ij}$ and uses it to define the number of scenarios $\sigma_i$ for each Subsystem $i$ in the main optimization. This result is used to formulate the locally discretized multi-subsystem staged-deployment optimization problem. This is similar to the fully flexible formulation, except that each subsystem evaluates a smaller number of scenarios through local scenario discretization, with the demand defined accordingly. See Appendix \ref{app} for more details on the demand definition with locally reduced scenarios. 

One detail that needs to be addressed in the formulation is the consistency constraint at the boundaries of each local scenario discretization. The last equation for each subsystem in Eq. \eqref{eq:full_example} takes the following form:
\begin{equation}
x_{ij}^{s} =x_{j}^{s}
\end{equation}
However, with the proposed local scenario discretization method, the number of scenarios considered for $x_{ij}$ and $x_j$ might not be the same anymore. Therefore, denoting the scenario set for $x_j$ by $\mathcal{S'}=\{1,2,3, ..., |S'|\}$ and the scenario set for $x_{ij}$ by $\mathcal{S''}=\{1,2,3, ..., |S''|\}$, the above constraint needs to be replaced by a consistency constraint $f$ in the following form:
\begin{equation}
x_{ij}^{s''} =f\left(x_{j}^{s'}\right), \quad s'\in \mathcal{S'},s''\in \mathcal{S''}
\end{equation}
In this step, we need to conserve the consistency of the scenarios so that, for example, in a special case where $|\mathcal{S''}|=|\mathcal{S'}|$, this relationship is reversible. Thus, in this mapping function $f$, we match the scenarios $s'$ and $s''$ as follows: 
\begin{equation} 
\label{consistency}
s''=\lceil{s'\cdot |\mathcal{S''}|/|\mathcal{S'}| }\rceil
\end{equation}

Consequently, the proposed locally discretized formulation reduces the number of weakly coupled variables and, in turn, the number of local scenarios at the subsystem level, while each subsystem still retains all the scenarios for its own demand uncertainties. Thus, we can maintain the flexibility in responding to the uncertain scenarios while limiting the computational effort. 

Note that this proposed approach includes multiple heuristic components developed for the staged-deployment problem like the example shown in Eq. \eqref{eq:full_example}. The statistical performance of the proposed approach, including all heuristic elements (e.g., the consistency mapping in Eq.~\eqref{consistency}), is demonstrated by the results in Sec. \ref{demo}. These results indicate that the collective effect of the heuristics does not introduce significant errors in the resulting solutions for the tested problems. The proposed approach may not be applicable to highly nonlinear or discrete problems, and a systematic study of the performance of this approach for more general stochastic programming problems is left for future work.

\subsection{Summary of the Method}
The following Algorithm \ref{alg1} summarizes the local scenario discretization method.

\begin{algorithm}
\caption{Local Scenario Discretization: 
}\label{alg1}
\begin{algorithmic}[1]

\State Run a deterministic optimization and obtain its solution $x^\text{det}$.
\State Compute the dynamic coupling $\delta_{ij}$ using the relationship between $\partial C_i/\partial D_{i}$ and $\partial C_i/\partial x_{ij}$ as well as $x^\text{det}$.
\State Assign the number of scenario $\sigma_i$ for each Subsystem $i$ using the dynamic coupling $\delta_{ij}$ and $x^\text{det}$.
\State Solve the staged-deployment optimization with the resulting formulation, where the consistency constraints apply when the number of scenarios changes at the subsystem boundaries.

\end{algorithmic}
\end{algorithm}

\section{Demonstration Case Studies and Analysis}
\label{demo}
This section demonstrates the developed local scenario discretization method for multi-subsystem staged co-deployment problems. We use the formulation in Eq. \eqref{eq:full_example} as our case studies and aim to compare the performance of the proposed method against the conventional methods under a fair setting. We first show a hypothetical example, then a computational experiment, and finally, a tutorial use case in the water-energy-food nexus application.

\subsection{Illustrative Example}
\label{demoexample}
We first consider a hypothetical example and perform a scalability analysis against the number of scenarios for each subsystem.  We use $c^1_A=4, c^2_A=5, c^1_B=1, c^2_B=2, d_{AB}=0.3, d_{BA}=0.1, \alpha=0.9, D^1_A=1, D^1_B=1, D^{2s}_A\sim U(1,2), D^{2s}_B\sim U(1,4)$, where $U(a,b)$ denotes a uniform distribution on an interval $[a,b]$. 
We set the number of scenarios for each Stage 2 demand $D^2_A$ and $D^2_B$ into $S$ scenarios (i.e., $S=S_A=S_B$), discretized with equal increments. (Note that defining the scenarios by random sampling from the probability distribution is also possible; however, that approach would lead to a high variance in the results when the number of scenarios is small \cite{kleywegt2002sample} and therefore is not reliable for large-scale optimization with limited computational resources.) With this definition of scenarios, the conventional fully flexible design would consider all $S_A \cdot S_B$ scenarios, whereas
the proposed method instead uses the developed local discretization method to reduce the computational effort leveraging the sensitivity information. 
The computational experiment is performed on a Windows 11 Enterprise machine with 12th Gen Intel(R) Core(TM) i7-1265U, 2.7 GHz, 10 cores, 12 logical processors. 
The formulation is converted into a nonlinear programming problem and solved using the \verb|fmincon| function in MATLAB R2022b, where the deterministic (i.e., inflexible and conservative) solution is used as an initial guess. The maximum number of function evaluations (\verb|MaxFunctionEvaluations|) in the \verb|fmincon| function is set as $200 \cdot N_\text{var}$, where $N_\text{var}$ is the number of variables for each optimization problem. The optimization is considered complete when \verb|exitflag| is nonnegative in \verb|fmincon|. Note that this test is conducted for a small-scale problem so that we can evaluate its performance by comparing it with the fully flexible solution; the proposed method can be applied to larger-scale problems where the fully flexible case is not computationally feasible.

Our goal is to evaluate the computational time and quality of the solution from the proposed method. The computational time is measured using the \verb|timeit| function in MATLAB. For the quality of the solution, note that the cost (i.e., the objective value from the optimization result) is not a fair metric to use; rather, we are interested in the Stage 1 decisions $x^1_A, x^1_B$ because they are the current decisions that we aim to optimize considering future uncertainties. One way to quantify the quality of our Stage 1 decisions against those from the fully flexible cases is to add a post-processing step that re-solves the fully flexible formulation while holding the Stage 1 variables fixed at the values from our proposed method (e.g., by setting both lower and upper bounds to the corresponding values). As a result, the optimality gap from the fully flexible cases over the same set of scenarios can be obtained for performance evaluation. Note that this step is only for evaluation purposes and is not used in practice and, therefore, is not included in computational time.

The comparison results between the deterministic (inflexible and conservative) case, the fully flexible case, and the proposed local scenario discretization method are shown in Table \ref{results}. From the results, there are a few clear trends. First, it is clear that the fully flexible case can significantly save costs compared to the deterministic case, reducing the (expected) cost significantly. In addition, the proposed local scenario discretization approach provides the Stage 1 decision results close to the fully flexible case, especially as the number of scenarios grows. This trend can be verified by examining the cost evaluated with fully flexible scenarios for the proposed method, which shows an optimality gap of less than 0.3\% of the fully flexible case (less than 0.1\% when the number of scenarios for each subsystem is four or more). Furthermore, the computational scalability of the proposed method is also evident, especially when the number of scenarios is large, because the proposed method reduces the problem size and significantly saves computational time. These results demonstrate the scalability of the proposed method, and the impact is expected to be more significant as we have more sources of uncertainty and scenarios.

\begin{table}[h!]

\centering
\normalsize\begin{tabulary}{1.0\textwidth}{CCCCC}\hline
Number of scenarios for each subsystem &  2 & 4 & 8 & 16 \\\hline
Deterministic (Conservative): Cost & \multicolumn{4}{c}{15.4510} \\
Deterministic (Conservative): $x^1_A$ & \multicolumn{4}{c}{3.2990} \\
Deterministic (Conservative): $x^1_B$ & \multicolumn{4}{c}{4.3300} \\
Deterministic (Conservative): Number of Variables & \multicolumn{4}{c}{ 6}\\
Deterministic (Conservative): Computational Time [s] & \multicolumn{4}{c}{ 0.0260}\\\hline
Fully Flexible: Cost & 13.5294 & 12.9712 & 12.7853 & 12.6979\\
Fully Flexible: $x^1_A$ & 2.2680 & 2.0275 & 2.1503 & 2.1518\\
Fully Flexible: $x^1_B$ &  4.2268 & 3.2303 & 2.9750 & 3.0319\\
Fully Flexible: Number of Variables & 18 & 66 & 258 & 1026\\
Fully Flexible: Computational Time [s] & 0.1080 & 0.2572 & 2.5784 & 106.5875
\\\hline
Proposed Method: Cost & 13.5683 & 13.0756 & 12.8430 & 12.7255\\
Proposed Method: $x^1_A$ &  2.2990 & 1.9895 & 2.1399 & 2.1485\\
Proposed Method: $x^1_B$ & 4.3299 & 3.1987 & 2.9750 & 3.0457\\
Proposed Method: Number of Variables &  14 & 34 & 114 & 418\\
Proposed Method: Computational Time [s]& 0.1006 & 0.1697 & 0.4598 & 17.6329\\
Proposed Method: Cost Evaluated with Fully Flexible Scenarios  & 13.5683 & 12.9793 & 12.7857 & 12.6926\\
Proposed Method: Optimality Gap from Fully Flexible Case  & 0.29\%& 0.06\%& 0.00\%& 0.04\%\\\hline
\end{tabulary}
\caption{{Computational results for the case study on the illustrative example.
{\label{results}}%
}}
\end{table}

\subsection{Computational Experiment over Random Problem Instances}
\label{experiment}
Due to the heuristics nature of the proposed method, it is of interest to evaluate its performance over a set of randomly generated problem instances. We consider a three-subsystem case and vary various parameters associated with the problem. Specifically, we consider a variety of problem instances with a linear cost model ($\alpha=1$) by sampling parameters $c^1_i, D^1_i, r^{c2}_i, r^{D2}_i$ from $U(0,1)$ so that $c^2_i=c^1_i\cdot (1+r^{c2}_i)$ and $D^{2s}_i \sim U(D^1_i,D^1_i\cdot (1+r^{D2}_i))$ for $i\in \{A,B,C\}$, as well as sampling $d_{ij}$ from $U(0,0.5)$ for $(i,j)\in \{A,B,C\}^2$. A linear programming (LP) tool Gurobi 11.0.1 is used for optimization. For the 1,000 sampled problem instances, we evaluate the optimality gap of the local scenario discretization method. 
The result of this computational experiment shows that the mean optimality gap is 0.0236\%, with a sample standard deviation of 0.0410\%, which shows the remarkable performance of the proposed local scenario discretization method.

\subsection{Use-Case Example: Water-Energy-Food Nexus}
\label{example}

In this subsection, we outline a use case for the methods described above in the application of the water-energy-food nexus. The purpose is to show a tutorial example of how the proposed method can be used. Rice farming is a motivating example for coupled subsystem design optimization under uncertainty. Macroscopically, rice production is dependent on available electricity for milling and water for rice paddies. Water is itself dependent on electricity because it must be pumped from the local water table. Furthermore, some forms of electricity generation require water. The design optimization of these systems cannot be done independently due to the coupling between them, and therefore is used as an example for coordinated staged deployment. Note that the numerical example in this subsection is presented as a use-case example; the parameter values and models are not validated.

We outline the parameters of the problem setup below. For this example analysis, we consider 8 scenarios for each subsystem. 
\begin{itemize}

\item We use labels A: Water generation system, assuming underground pumping; B: Electricity generation system, assuming thermoelectric power plant; C: Rice production system.

\item 1 unit of demand for each subsystem is defined using units of 1 m$^3$ of water, 1 kWh of electricity, and 1 kg of rice.

\item For the cost model of infrastructure deployment, we use the cost of resources as proxies. A linear cost model is considered ($\alpha=1$).

\item Costs for Stage 1: $c^1_A=2.16$ USD per unit; $c^1_B=0.163$ USD per unit; $c^1_C=2$ USD per unit.

\item Costs for Stage 2: $c^2_A, c^2_B, c^2_C$ are $1.06^{15}$ times $c^1_A, c^1_B, c^1_C$, assuming a 6\% increase per year for 15 years, consistent with the current global inflation rate as of 2023.

\item Demands for Stage 1: $D^1_A=1385$ units, per capita per year; $D^1_B=2940$ units, per capita per year; $D^1_C=80$ units, per capita per year.

\item Demands for Stage 2: $D^{2s}_A\sim U(D^1_A,1.2D^1_A), D^{2s}_B\sim U(D^1_B,2D^1_B) , D^{2s}_C\sim  U(D^1_C,1.4D^1_C)$.

\item Coupling between subsystems: $d_{AB}=0.057$; $d_{BA}=0.147$; $d_{AC}=2.5$; $d_{CA}=1 \times 10^{-10}$; $d_{BC}=0.05$;  $d_{CB}=1 \times 10^{-10}$. The value $1 \times 10^{-10}$ is used to represent the smallest amount of coupling. This means: to generate 1 unit of water ($1 \text{m}^3$), one consumes 0.147 units ($0.147 \  \text{kWh}$) of electricity; to generate 1 unit of electricity ($1 \ \text{kWh}$) one consumes 0.057 units ($0.057 \  \text{m}^3$) of water (e.g., due to evaporation); and to generate 1 unit of rice ($1 \text{kg}$), one consumes 2.5 units ($2.5 \text{m}^3$)  of water and 0.05 units ($0.05\  \text{kWh}$) of electricity.
\end{itemize}





Table \ref{tab:result6} shows the results for this use case. Note that, for this particular problem, Subsystem C can be solved separately due to the negligible values of $d_{CA}$ and $d_{CB}$, but the results shown here are from the integrated optimization for demonstration purposes.
From the results, it is clear that the local scenario discretization method shows an excellent approximation to the original problem. While we do not claim that this numerical case setup or results are representative of a real-world case, we hope that this use-case tutorial provides users with a guideline on how to apply the proposed framework to their use cases.

\begin{table}[h!]

\centering
\normalsize\begin{tabulary}{1.0\textwidth}{CCC}\hline
Variable & Fully Flexible & Local Scenario Discretization  \\\hline
 $x^1_A$ & 2058.4261 & 2061.8483
\\
$x^1_B$ & 4933.1405 & 4923.3329
\\
$x^1_C$ & 98.2857 & 98.2857
\\\hline
Optimality Gap&-- &  0.0026\% 
\\\hline
\end{tabulary}
\caption{{Computational results for the water-energy-food nexus application. 
{\label{tab:result6}}%
}}
\end{table}

\section{Extensions}
\label{discussion}
The developed method provides a new perspective to capture the coupling between subsystems in a design optimization problem. This section discusses the implications of our method in the context of system partitioning for optimization of complex system design.

\subsection{System Partitioning via Dynamic Coupling}
The flexible staged-deployment design is challenging due to the coupled subsystems; the examples tested in the previous section can be solved in a monolithic formulation, but that formulation can become prohibitively expensive to solve when there are a larger number of subsystems. One approach to solving such complex and intractable problems is to partition systems into subsystem clusters and solve them using a decomposition-based optimization approach \cite{yi2008comparison,agte2010mdo,martins2013multidisciplinary}. 

System partition optimization for decomposition-based optimization generally involves a trade-off between the complexity of partitioned sub-problems and the coordination effort \cite{allison2009optimal}. Mathematically, this problem can be formulated as a multi-objective optimization problem: $\min_{\boldsymbol{p}}(CS,SS_\text{max})$, where $\boldsymbol{p}$ is the partitioning vector (i.e., which subsystem belongs to which subproblem), $CS$ quantifies the coordination effort, and $SS_\text{max}$ is the maximum subproblem size. In the past literature, $CS$ was defined using the number of shared and coupling variables, and $SS_i$ for Subsystem $i$ was quantified using the number of associated decision variables, consistency constraints, and/or analysis functions. In practice, we often expect to minimize $CS$ within the computational resource constraints, and therefore a more practical formulation could be: $\min_{\boldsymbol{p}}CS \quad \text{subject to} \quad SS_\text{max}\leq SS^\text{UB}$, where $SS^\text{UB}$ is the maximum size of the subproblem given the computational resource constraints.
A caveat of this approach, when applied to our system co-deployment problem, is that it does not take into consideration the dynamic coupling between the subsystems.

Our local scenario discretization provides a new way to perform the system partitioning for the system co-deployment problem. Namely, we can use the dynamic coupling metric, a byproduct of our local scenario discretization method, to capture the coordination complexity more explicitly and introduce a new definition for $CS$ and $SS_\text{max}$. For example, for a pair of subsystems that have one coupling variable $x_{ij}$, the conventional system partitioning approach would just count that as one coupling variable in computing $CS$, whereas our approach can actually evaluate the dynamic coupling strength $\delta_{ij}$ and use it to represent $CS$, more accurately characterizing the coordination effort. Thus, one example definition is to define $CS$ as the summation of $\delta_{ij}$ across the partitioning boundary and $SS_\text{max}$ as the maximum number of variables in a (locally-scenario-discretized) subproblem. In this way, we can consider both the extent of couplings and their temporal evolution in partitioning; this new partitioning problem is suitable for the flexible design problem because the coordination effort for the highly scenario-dependent coupling variables would be heavier than that for scenario-independent coupling variables.

\subsection{Bottom-Up Built-In Flexibility Design with the Partitioned Subproblems}
Beyond system partitioning, our local scenario discretization also enables a system designer to design and evaluate built-in flexibility from the bottom up. Conventional staged-deployment design approaches often analyze and identify the sources of flexibility at a high level and then analyze the impact of these identified sources with more detailed models \cite{wang2005real}, because optimizing all decisions capturing all uncertainty scenarios can be intractable or impractical for complex systems. In contrast, the proposed approach, especially when used with system partitioning, enables us to optimize staged-deployment decisions within a partitioned subproblem with strongly-connected subsystems only. This partitioning can be handled by representing the couplings from outside of the partitioned subproblem using the deterministic (conservative) solution, which can be obtained efficiently prior to the main optimization, while solving for the uncertainties within the subproblem in a stochastic way. 

As a demonstrating example, consider an extension of the illustrative example in Section \ref{demoexample} to a three-subsystems case with parameters: $c^1_A=4, c^2_A=5, c^1_B=1, c^2_B=2, c^1_C=2, c^2_C=3, d_{AB}=0.3, d_{BA}=0.1, d_{AC}=0.03, d_{CA}=0.03, d_{BC}=0.01, d_{CB}=0.01, \alpha=0.9, D^1_A=1, D^1_B=1, D^1_C=1, D^{2s}_A\sim U(1,2), D^{2s}_B\sim U(1,4), D^{2s}_C\sim U(1,2)$. We consider 8 scenarios for each subsystem in this example. For this problem, we consider four partitioning options: AB-C, BC-A, AC-B, and A-B-C. Here, XY-Z, for example, indicates that Subsystems X and Y are partitioned into one subproblem and Subsystem Z into the other subproblem. In this problem, for the deterministic case, both $SS_{max}$ and $CS$ are the same for all three partitioning options. However, when the local scenario discretization is considered, $CS$, defined as the summation of $\delta_{ij}$ across the boundary of the partition, and $SS_{max}$, defined as the maximum number of variables in the subproblem, would be different for each partitioning option. 

Table \ref{results2} shows the results for different partitioning options, including their $CS$ and $SS_\text{max}$ as well as their Stage 1 decision variables and optimality gap. 
We can see that $CS$ is smaller for AB-C than for the other partitioning options due to the weak coupling between Subsystem C and the rest of the system, whereas $SS_{max}$ is largest for AB-C. 
Assuming $SS_\text{max}$ for any of these options is within the computational constraints, we can pick AB-C as our optimal partitioning.
The result shows that, with this partitioning option AB-C, the bottom-up approach gives Stage 1 decisions that are close to the optimal solution given by the fully flexible case and thus with a small optimality gap.
The other partitioning options BC-A, AC-B, and A-B-C are not as effective as the optimal partitioning AB-C. 

\begin{table}[h!]

\centering
\normalsize\begin{tabulary}{1.0\textwidth}{CCCCCCC}\hline
 Variable & Fully Flexible & Local Scenario Discretization  & AB-C & BC-A & AC-B & A-B-C\\\hline
 $CS$ & -- & -- & 0.1099 & 0.8208 & 0.7653 & 0.8480\\
 $SS_\text{max}$ & -- & -- & 114 & 34 & 34 & 17\\\hline
$x^1_A$ & 2.2110 & 2.2159 & 2.2129 & 2.8005 & 2.7828 & 2.8005\\
$x^1_B$ & 3.0070 & 2.9989 & 3.0037 & 3.0685 & 3.0729 & 3.0729\\
$x^1_C$ &  1.5496 & 1.5404 & 1.5733 & 1.5561 & 1.5557 & 1.5733
\\\hline
Optimality Gap &   -- &  0.42\% &  0.39\% &  3.85\% &  3.69\% &  3.87\% \\\hline
\end{tabulary}
\caption{{Computational results for all problem partitioning options. 
{\label{results2}}%
}}
\end{table}

This bottom-up flexibility design method is not necessarily universally effective, but rather is particularly suitable when one or more of the subsystems can be effectively decoupled from the rest of the subsystem. To show this, we further perform two sets of computational experiments over random problem instances: Case 1: all subsystems equally coupled, using the same parameter values and range as in Section \ref{experiment}; and Case 2: Subsystem C is relatively decoupled from the rest, constructed by reducing $d_{ij}$ by a factor of 10 compared to Case 1 when $i \ \text{or} \ j=C$.\footnote{Note that a reduced value of $d_{ij}$ when $i \ \text{or} \ j=C$ does not automatically make AB-C the optimal partitioning option; while AB-C has a higher chance of being optimal, the actual partitioning performance also depends on the rest of the problem.}
With 1,000 uniformly sampled problem instances for each of Cases 1 and 2, the mean optimality gap with the lowest-$CS$ partition is 2.2013\% for Case 1 and 0.1275\% for Case 2.
These optimality gaps are higher than those with the local scenario discretization method (0.0236\% for Case 1 and 0.0261\% for Case 2) due to the lost information in the partitioning process, but the advantage of the bottom-up flexibility design method is that it does not need to solve the integrated problem (which can be computationally challenging or impractical) while still achieving an acceptable level of optimality gap, especially for Case 2.
Note that, due to its heuristic nature, the proposed partitioning method based on the $CS$ values does not necessarily select the actual lowest optimality gap of all partitioning options; in fact, the probability that the proposed method selects the actual lowest optimality gap within a tolerance of 0.01\% (chosen as a reasonably acceptable level in practice) is 58.9\% of Case 1 and 90.5\% for Case 2.
This trend is expected, as it validates our hypotheses that the proposed bottom-up built-in flexibility design works well for Case 2 (where a subsystem is decoupled from the rest), but not necessarily for Case 1 (where all subsystems are equally well-connected). 
Even when the proposed method does not find the optimal (or even feasible) solution for the original problem, the found solutions (even the infeasible ones) can still be used as good initial solutions for further refinement using decomposition-based optimization techniques \cite{yi2008comparison,agte2010mdo,martins2013multidisciplinary}. 

Finally, it is worth noting that, although the above example only considers the Stage 1 decisions as parameters to be optimized from the bottom up, the same approach can be applied to any parameters that govern the subsystem decisions in response to uncertainties (e.g., parameters in decision rules). This bottom-up approach expands the design-for-flexibility method to more complex systems than the conventional approaches, providing effective solutions even when the conventional approaches are challenging or intractable.

\section{Conclusion}
\label{conclusion}
In this paper, the nonlinear stochastic programming formulation for staged-deployment design is extended to complex multi-subsystem systems, and a new local scenario discretization method is developed. While the conventional stochastic programming method for staged deployment uses the same scenario discretization across the entire system, which is not scalable, the developed local scenario discretization method uses a different number of scenarios for each subsystem by effectively decoupling weak dynamic coupling between subsystems. This is enabled by characterizing the strength of the dynamic coupling between subsystems, which is then used to determine the number of scenarios considered for each coupling variable. 
The resulting formulation and method are tested with case studies, and their efficient computational performance is demonstrated. Extensions of the proposed method in the context of system partitioning are also discussed.
This method is expected to improve the staged-deployment design of various complex engineering systems under uncertainty, such as water, energy, food, and other infrastructure systems.

\section*{Acknowledgments}\label{acknowledgements}
This material is based upon work supported by the National Science Foundation under Grant: 1942559. 
AI technologies (Microsoft Copilot, with occasional use of ChatGPT and Google Gemini) were used for technical and editorial support. Specifically, these tools assisted with suggesting candidate parameter ranges for the use-case example, drafting auxiliary interface code, supporting summaries of relevant literature for review, and improving language clarity. 
Any AI-assisted text or code was reviewed, edited, and validated by the authors, who take full responsibility for the final content.
AI technologies did not generate original scientific findings or replace the authors’ critical analysis.
The authors would like to thank Yuhei Miyauchi for discussions on an earlier version of the method. The authors would also like to thank Zachary Grieser and Nicholas Gollins for their support in reviewing the code and manuscript.



\appendix
\section{Appendix A: Demand Definition for Scenarios}
\label{app}
This Appendix discusses the definition of the demand for the fully flexible and proposed locally discretized formulations. There are a variety of different ways to define and implement the demand scenarios, but the following shows one way. 

Consider a system with $N$ subsystems. For simplicity, we consider a case where each subsystem has one uncertain parameter, and assume that each subsystem's uncertain parameter (e.g., demand) has $S_i$ scenarios. 

In the fully flexible case, all variables need to consider all combinations of scenarios. For Subsystem $i$, the demand has $S_i$ scenarios ($D^1, D^2, D^3, ..., D^{S_i}$), but its decision variables need to consider the entire $\prod_j{S_j}$ scenarios. Thus, we define the demand for each scenario in the following way: each entry of the demand repeats $\prod_{j=1}^{i-1}S_j$ times to form a pattern, and that pattern repeats $\prod_{j=i+1}^{N}S_j$ times. For example, if there are $3$ subsystems and each has its demand discretized into two scenarios of $D=1,2$. In this case, there would be $2 \cdot 2 \cdot 2 = 8$ scenarios for all variables. Using the above approach, the demand of subsystem 1 for these $8$ scenarios would be: $1,2,1,2,1,2,1,2$, respectively, that of Subsystem 2 would be $1,1,2,2,1,1,2,2$, respectively, and that of Subsystem 3 would be $1,1,1,1,2,2,2,2$, respectively. This definition allows us to systematically explore the entire multi-dimensional uncertainty space. The mathematical equation can be coded using \verb|repmat| and \verb|repelem| functions (or \verb|combvec| function) in MATLAB. 

The above definition can naturally be extended to the locally discretized formulation. Here, the number of scenarios differs for each system's decision variables. Specifically, the decision variables in Subsystem $i$ consider $S_i\cdot\prod_{j:j\neq i}{S_{ij}}$ scenarios, where $S_{ij}$ is the contribution by the coupling variables from subsystems $j$ to $i$. Thus, we consider the following patterns: each entry of the demand repeats $\prod_{j=1}^{i-1}S_{ij}$ times to form a pattern, and that pattern repeats $\prod_{j=i+1}^{N}S_{ij}$ times. Similarly to the fully flexible case, this pattern can be formulated using \verb|repmat| and \verb|repelem| functions in MATLAB.

\bibliography{bibliography/converted_to_latex.bib%
}

@book{neufville2011flexibility,
  title={Flexibility in Engineering Design. Engineering Systems},
  author={Neufville, Richard de and Scholtes, Stefan},
  year={2011},
  publisher={Cambridge, MA, United States: MIT Press}
}

@article{de2004staged,
  title={Staged deployment of communications satellite constellations in low earth orbit},
  author={de Weck, Olivier L and De Neufville, Richard and Chaize, Mathieu},
  journal={Journal of Aerospace Computing, Information, and Communication},
  volume={1},
  number={3},
  pages={119--136},
  year={2004},
  doi={10.2514/1.6346}
}

@article{hamdan2023dynamic,
  title={A dynamic multi-stage design framework for staged deployment optimization of highly stochastic systems},
  author={Hamdan, Bayan and Liu, Zheng and Ho, Koki and B{\"u}y{\"u}ktahtak{\i}n, {\.I} Esra and Wang, Pingfeng},
  journal={Structural and Multidisciplinary Optimization},
  volume={66},
  number={7},
  pages={162},
  year={2023},
  publisher={Springer},
  doi={10.1007/s00158-023-03609-6}
}

@article{caunhye2017approach,
  title={An approach based on robust optimization and decision rules for analyzing real options in engineering systems design},
  author={Caunhye, Aakil M and Cardin, Michel-Alexandre},
  journal={IISE Transactions},
  volume={49},
  number={8},
  pages={753--767},
  year={2017},
  publisher={Taylor \& Francis},
  doi={10.1080/24725854.2017.1299958}
}

@article{cardin2017approach,
  title={An approach for analyzing and managing flexibility in engineering systems design based on decision rules and multistage stochastic programming},
  author={Cardin, Michel-Alexandre and Xie, Qihui and Ng, Tsan Sheng and Wang, Shuming and Hu, Junfei},
  journal={IISE Transactions},
  volume={49},
  number={1},
  pages={1--12},
  year={2017},
  publisher={Taylor \& Francis},
  doi={10.1080/0740817X.2016.1189627}
}

@inproceedings{hassan2005value,
  title={Value-at-risk analysis for real options in complex engineered systems},
  author={Hassan, Rania and De Neufville, Richard and McKinnon, Douglas},
  booktitle={2005 IEEE International Conference on Systems, Man and Cybernetics},
  volume={4},
  pages={3697--3704},
  year={2005},
  organization={IEEE},
  doi={10.1109/ICSMC.2005.1571721}
}

@article{engel2008designing,
  title={Designing systems for adaptability by means of architecture options},
  author={Engel, Avner and Browning, Tyson R},
  journal={Systems Engineering},
  volume={11},
  number={2},
  pages={125--146},
  year={2008},
  publisher={Wiley Online Library},
  doi={10.1002/sys.20090}
}

@article{siddiqi2011water,
  title={The water--energy nexus in Middle East and North Africa},
  author={Siddiqi, Afreen and Anadon, Laura Diaz},
  journal={Energy Policy},
  volume={39},
  number={8},
  pages={4529--4540},
  year={2011},
  publisher={Elsevier},
  doi={10.1016/j.enpol.2011.04.023}
}

@book{powell2007approximate,
  title={Approximate Dynamic Programming: Solving the curses of dimensionality},
  author={Powell, Warren B},
  year={2011},
  publisher={John Wiley \& Sons},
  doi={10.1002/9781118029176}
}

@article{dupavcova2003scenario,
  title={Scenario reduction in stochastic programming},
  author={Dupa{\v{c}}ov{\'a}, Jitka and Gr{\"o}we-Kuska, Nicole and R{\"o}misch, Werner},
  journal={Mathematical Programming},
  volume={95},
  pages={493--511},
  year={2003},
  publisher={Springer},
  doi={10.1007/s10107-002-0331-0}
}

@article{kleywegt2002sample,
  title={The sample average approximation method for stochastic discrete optimization},
  author={Kleywegt, Anton J and Shapiro, Alexander and Homem-de-Mello, Tito},
  journal={SIAM Journal on Optimization},
  volume={12},
  number={2},
  pages={479--502},
  year={2002},
  publisher={SIAM},
  doi={10.1137/S1052623499363220}
}

@phdthesis{wang2005real,
  title={Real options "in" projects and systems design: identification of options and solutions for path dependency},
  author={Wang, Tao},
  year={2005},
  school={Massachusetts Institute of Technology}
}

@article{allison2009optimal,
    author = {Allison, James T. and Kokkolaras, Michael and Papalambros, Panos Y.},
    title = {Optimal Partitioning and Coordination Decisions in Decomposition-Based Design Optimization},
    journal = {Journal of Mechanical Design},
    volume = {131},
    number = {8},
    pages = {081008},
    year = {2009},
    month = {07},
    doi = {10.1115/1.3178729},
    eprint = {https://asmedigitalcollection.asme.org/mechanicaldesign/article-pdf/131/8/081008/5524404/081008\_1.pdf},
}

@article{agte2010mdo,
  title={MDO: assessment and direction for advancement—an opinion of one international group},
  author={Agte, Jeremy and De Weck, Olivier and Sobieszczanski-Sobieski, Jaroslaw and Arendsen, Paul and Morris, Alan and Spieck, Martin},
  journal={Structural and Multidisciplinary Optimization},
  volume={40},
  pages={17--33},
  year={2010},
  publisher={Springer},
  doi={10.1007/s00158-009-0381-5}
}

@article{homem2014monte,
  title={Monte Carlo sampling-based methods for stochastic optimization},
  author={Homem-de-Mello, Tito and Bayraksan, G{\"u}zin},
  journal={Surveys in Operations Research and Management Science},
  volume={19},
  number={1},
  pages={56--85},
  year={2014},
  publisher={Elsevier},
  doi={10.1016/j.sorms.2014.05.001}
}

@incollection{homem2014stochastic,
  title={Stochastic constraints and variance reduction techniques},
  author={Homem-de-Mello, Tito and Bayraksan, G{\"u}zin},
  booktitle={Handbook of Simulation Optimization},
  pages={245--276},
  year={2014},
  publisher={Springer},
  doi={10.1007/978-1-4939-1384-8_9}
}

@article{martins2013multidisciplinary,
  title={Multidisciplinary design optimization: a survey of architectures},
  author={Martins, Joaquim RRA and Lambe, Andrew B},
  journal={AIAA journal},
  volume={51},
  number={9},
  pages={2049--2075},
  year={2013},
  publisher={American Institute of Aeronautics and Astronautics},
  doi={10.2514/1.J051895}
}

@book{shapiro2021lec_on_SP,
author = {Shapiro, Alexander and Dentcheva, Darinka and Ruszczynski, Andrzej},
title = {Lectures on Stochastic Programming: Modeling and Theory, Third Edition},
publisher = {Society for Industrial and Applied Mathematics},
year = {2021},
doi = {10.1137/1.9781611976595},
address = {Philadelphia, PA},
edition   = {3},
}

@article{ahmed2003MultiStageStochasticInteger,
  title = {A Multi-Stage Stochastic Integer Programming Approach for Capacity Expansion Under Uncertainty},
  author = {Ahmed, Shabbir and King, Alan J. and Parija, Gyana},
  year = {2003},
  journal = {Journal of Global Optimization},
  volume = {26},
  number = {1},
  pages = {3--24},
  doi = {10.1023/A:1023062915106},
}

@article{ahmed2004TSSIP_BnB,
  title = {A Finite Branch-and-Bound Algorithm for Two-Stage Stochastic Integer Programs},
  author = {Ahmed, Shabbir and Tawarmalani, Mohit and Sahinidis, Nikolaos V.},
  year = {2004},
  month = jun,
  journal = {Mathematical Programming},
  volume = {100},
  number = {2},
  pages = {355--377},
  doi = {10.1007/s10107-003-0475-6},
}

@article{vanderlaan2023MIRbenders,
    author = {van der Laan, Niels and Romeijnders, Ward},
    title = {A Converging Benders’ Decomposition Algorithm for Two-Stage Mixed-Integer Recourse Models},
    journal = {Operations Research},
    volume = {72},
    number = {5},
    pages = {2190--2214},
    year = {2024},
    doi = {10.1287/opre.2021.2223},
}

@article{yi2008comparison,
  title={Comparison of MDO methods with mathematical examples},
  author={Yi, Sang-Il and Shin, Jung-Kyu and Park, GJ},
  journal={Structural and Multidisciplinary Optimization},
  volume={35},
  number={5},
  pages={391--402},
  year={2008},
  publisher={Springer},
  doi={10.1007/s00158-007-0150-2}
}

@article{bertsimas2023optimization,
  title={Optimization-based scenario reduction for data-driven two-stage stochastic optimization},
  author={Bertsimas, Dimitris and Mundru, Nishanth},
  journal={Operations Research},
  volume={71},
  number={4},
  pages={1343--1361},
  year={2023},
  publisher={INFORMS},
  doi={10.1287/opre.2022.2265}
}

@article{henrion2022problem,
  title={Problem-based optimal scenario generation and reduction in stochastic programming},
  author={Henrion, R{\'e}ne and R{\"o}misch, Werner},
  journal={Mathematical Programming},
  volume={191},
  number={1},
  pages={183--205},
  year={2022},
  publisher={Springer},
  doi={10.1007/s10107-018-1337-6}
}

@article{hewitt2022decision,
  title={Decision-based scenario clustering for decision-making under uncertainty},
  author={Hewitt, Mike and Ortmann, Janosch and Rei, Walter},
  journal={Annals of Operations Research},
  volume={315},
  number={2},
  pages={747--771},
  year={2022},
  publisher={Springer},
  doi={10.1007/s10479-020-03843-x}
}

@article{heitsch2003scenario,
  title={Scenario reduction algorithms in stochastic programming},
  author={Heitsch, Holger and R{\"o}misch, Werner},
  journal={Computational Optimization and Applications},
  volume={24},
  number={2},
  pages={187--206},
  year={2003},
  publisher={Springer},
  doi={10.1023/A:1021805924152}
}

@book{pflug2014multistage,
  title={Multistage Stochastic Optimization},
  author={Pflug, Georg Ch and Pichler, Alois},
  year={2014},
  publisher={Springer},
  doi={10.1007/978-3-319-08843-3}
}

@article{rujeerapaiboon2022scenario,
  title={Scenario reduction revisited: Fundamental limits and guarantees},
  author={Rujeerapaiboon, Napat and Schindler, Kilian and Kuhn, Daniel and Wiesemann, Wolfram},
  journal={Mathematical Programming},
  volume={191},
  number={1},
  pages={207--242},
  year={2022},
  publisher={Springer},
  doi={10.1007/s10107-018-1269-1}
}

@article{kim2014guide,
  title={A guide to sample average approximation},
  author={Kim, Sujin and Pasupathy, Raghu and Henderson, Shane G},
  journal={Handbook of simulation optimization},
  pages={207--243},
  year={2014},
  publisher={Springer},
  doi={10.1007/978-1-4939-1384-8_8}
}

@article{hoyland2003heuristic,
  title={A heuristic for moment-matching scenario generation},
  author={H{\o}yland, Kjetil and Kaut, Michal and Wallace, Stein W},
  journal={Computational optimization and applications},
  volume={24},
  number={2},
  pages={169--185},
  year={2003},
  publisher={Springer},
  doi={10.1023/A:1021853807313}
}

@article{suemitsu2025surrogate,
  title={Surrogate-Assisted Scenario-Generation Method for Simulation-Based Stochastic Programming Problems},
  author={Suemitsu, Issei and Izui, Kazuhiro},
  journal={IEEE Transactions on Automation Science and Engineering},
  volume={22},
  pages={13161--13174},
  year={2025},
  publisher={IEEE},
  doi={10.1109/TASE.2025.3551716}
}

@article{patel2022neur2sp,
  title={Neur2SP: Neural two-stage stochastic programming},
  author={Dumouchelle, Justin and Patel, Rahul Mihir and Khalil, Elias and Bodur, Merve},
  journal={Advances in neural information processing systems},
  volume={35},
  pages={23992--24005},
  year={2022}
}

@inproceedings{wu2022learning,
  title={Learning scenario representation for solving two-stage stochastic integer programs},
  author={Wu, Yaoxin and Song, Wen and Cao, Zhiguang and Zhang, Jie},
  year={2022},
  organization={International Conference on Learning Representations, ICLR}
}

\end{document}